\documentclass{article}
\begin{document}
\begin{center}\Large{Counting on Rectangular Areas}\end{center}
 \begin{center}Milan Janji\'c,\end{center}
 \begin{center}Faculty of Natural Sciences and mathematics, \end{center}
\begin{center} Banja Luka, Republic of Srpska,
Bosnia and Herzegovina.\end{center}

\begin{center}\Large{Counting on Rectangular Areas}\end{center}

\begin{abstract}
In the first section of this paper we prove a theorem  for the
number of columns of a rectangular area that are identical to the
given one. A  special case, concerning $(0,1)$-matrices, is also
stated.

In the next section we apply this theorem to  derive several
combinatorial identities
 by counting specified subsets of  a finite set. This means that the obtained identities
 will involve
  binomial coefficients only.
We start with a simple equation
 which is, in fact, an immediate
consequence of Binomial theorem, but it is derived independently
of it. The second result concerns sums of binomial coefficients.
In a special case we obtain one of the best known binomial
identity
 dealing with alternating sums. Klee's identity is also obtained as a special case as well as some formulae
 for partial sums of binomial coefficients, that is, for the numbers of
 Bernoulli's triangle.

\end{abstract}
\section{A counting theorem}

The set of natural numbers $\{1,2,\ldots,n\}$ will be denoted by
$[n],$ and by $|X|$ will be denoted the number of elements of the
set $X.$

For the proof of the main theorem we need the following simple
result:
\begin{equation}\label{u1}\sum_{I}(-1)^{|I|}=0,\end{equation}
where $I$ run over all subsets of $[n]$ (empty set included). This
may be easily proved by induction or using Binomial theorem. But
the proof by induction makes all further investigations
independent even of Binomial theorem.

 Let $A$ be an
$m\times n$ rectangular matrix filled with elements which belong
to a set $\Omega.$

 By the i-column of $A$ we shall mean each column of $A$ that is equal
 to $[c_1,c_2,\ldots,c_m]^T,$
where  $c_1,c_2,\ldots,c_m$ of $\Omega$ are given. We shall denote
the number of i-columns of $A$  by $\nu_A(c)$ or simply by
$\nu(c).$

For $I=\{i_1,i_2,\ldots,i_k\}\subset [m],$
 by $A(I)$ will be denoted the maximal number of columns $j$ of $A$  such that
$$a_{ij}\not=c_j,\;(i\in I).$$

We also define $$A(\emptyset)=n.$$

\textbf{Theorem 1.}\label{t1}\textit{ The number $\nu(c)$ of
i-columns of $A$ is equal
\begin{equation}\label{tt1}\nu(c)=\sum_{I}(-1)^{|I|}A(I),\end{equation}
 where summation is taken over all  subsets $I$ of
$[m].$}

 Proof. Theorem may be proved by the standard combinatorial method, by counting
 the contribution of each column of $A$ in the sum
 on the right side of (\ref{tt1}).

We give here a proof by induction. First, the formula will be
proved in the case $\nu(c)=0$ and $\nu(c)=n.$ In the case
$\nu(c)=n$ it is obvious that for $I\not=\emptyset$ we have
$A(I)=0,$ which implies
$$\sum_{I}(-1)^{|I|}A(I)=n+\sum_{I\not=\emptyset}(-1)^{|I|}A(I)=n.$$

In the case $\nu(c)=0$ we use induction on $n.$ If $n=1$ then the
matrix $A$ has only one  column, which is not equal $c.$ It yields
that there exists $i_0\in\{1,2,\ldots,m\}$ such that
$a_{i_0,1}\not=c_{i_0}.$
 Denote by $I_0$ the set of all such numbers. Then   $A(I)=1$ if and only if $I\subset I_0.$
 From this and (\ref{u1})
 we obtain
$$\sum_{I}(-1)^{|I|}A(I)=\sum_{I\subset I_0}(-1)^{|I|}=0.$$

Suppose now that the formula is true for matrices with $n$ columns
and that $A$ has $n+1$-columns, and $\nu_A(c)=0.$  Omitting the
first column, the matrix $B$ with $n$ columns remains. If $I_0$ is
the same as in the case $n=1,$ then
$$\sum_{I}(-1)^{|I|}A(I)=\sum_{I\not\subset I_0}(-1)^{|I|}A(I)+
\sum_{I\subset I_0}(-1)^{|I|}A(I)=$$$$=\sum_{I\not\subset
I_0}(-1)^{|I|}B(I) +\sum_{I\subset
I_0}(-1)^{|I|}(B(I)+1)=$$$$=\sum_{I}(-1)^{|I|}B(I)+\sum_{I\subset
I_0}(-1)^{|I|}=0,$$ since the first sum is equal zero by the
induction hypothesis, and the second by (\ref{u1}).

For the rest of the proof we use induction on $n$ again. For $n=1$
the matrix $A$ has only one column which is either equal  $c$ or
not. In both cases theorem  is true, from the preceding.

 Suppose
that theorem holds for $n,$ and that the matrix $A$ has $n+1$
columns. We may suppose that $\nu(c)\geq 1.$  Omitting one of the
i-columns  we obtain the matrix $B$ with $n$ columns.
 By the induction hypothesis theorem is true for $B$.
On the other hand it is clear that $A(I)=B(I)$ for each nonempty
subset $I.$ Furthermore $A$ has one i-column more then $B,$ which
implies

$$\nu(c)=\nu_A(c)=\nu_B(c)+1=1+\sum_I(-1)^{|I|}B(I)=$$$$=1+n+\sum_{I\not=\emptyset}(-1)^{|I|}B(I)
=1+n+\sum_{I\not=\emptyset}(-1)^{|I|}A(I).$$ Thus
$$\nu(c)=\sum_{I}(-1)^{|I|}A(I),$$ and theorem is proved.

If the number $A(I)$ does not depend on elements of the set $I,$
but only on its number $|I|$
 then the equation(\ref{tt1}) may be written in the form
\begin{equation}\label{ni}\nu(c)=\sum_{i=0}^m(-1)^i{m\choose i}A(i),\end{equation}
where $|I|=i.$

Our object of investigation will be $(0,1)$ matrices.  Let $c$ be
the i- column of a such matrix  $A.$ Take $I_0\subseteq
[m],\;|I_0|=k$ such that
\begin{equation}\label{ck}c_i=\left\{\begin{array}{cc}1&i\in
I_0\\0&i\not\in I_0\end{array}\right.\end{equation}

Then the number $A(I)$ is equal to the number of columns of $A$
having  $0$'s in the rows labelled by the set $I\cap I_0,$ and
$1$'s in the rows labelled by the set $I\setminus I_0.$ Suppose
that the number $A(I)$ depends only on
 $|I\cap I_0|,\;|I\setminus I_0|.$ If we denote  $|I\cap I_0|=i_1,\;|I\setminus I_0|=i_2,\;A(I)=A(i_1,i_2),$
then (\ref{tt1}) may be written in the form
\begin{equation}\label{rv}\nu(c)=\sum_{i_1=0}^{k}\sum_{i_2=0}^{m-k}(-1)^{i_1+i_2}{k\choose
i_1}{m-k\choose i_2}A(i_1,i_2).\end{equation}

\section{Counting subsets of a finite set}
Suppose that a finite set $X=\{x_1,x_2,\ldots, x_n\}$ is given.
Label by $1,2,\ldots, 2^n$ all subsets of $X$ arbitrary and define
an $n\times 2^n$ matrix $A$ in the following way
\begin{equation}\label{m1}a_{ij}=\left\{\begin{array}{cc}1&\mbox{ if $x_i$
lies in the set labelled by }j\\0&\mbox{ otherwise
}\end{array}\right..\end{equation}

Take $I_0\subseteq [n],\;|I_0|=k,$ and form the submatrix $B$ of
$A$ consisting of those rows of $A$ which indices belong to $I_0.$
Let $c$ be arbitrary i-column of $B.$ Define
\begin{equation}\begin{array}{c}I_0'=\{i\in I_0:c_i=1\},\\I_0''=\{i\in I_0:c_i=0\}\end{array}.\end{equation}

The number $\nu(c)$ is equal to the number of subsets that contain
the set $\{x_i,\;i\in I_0'\},$ and do not intersect the set
$\{x_i:i\in I_0''\}.$ There are obviously $$\nu(c)=2^{n-k},$$ such
sets.

Furthermore, if $I\subseteq I_0$ then the number $B(I)$ is equal
to the number of subsets that contain the set $\{x_i:i\in I\cap
I_0''\},$ and do not meet the set $\{x_i:i\in I\cap I_0'\}.$ It is
clear that there are  $$B(I)=2^{n-|I|}$$ such subsets, so that the
formula (\ref{tt1}) may be applied. It follows
$$2^{n-k}=\sum_{i=0}^k(-1)^i{k\choose i}2^{n-i}.$$

Thus we have

\noindent\textbf{Proposition 2.1.} \textit{For each nonnegative
integer $k$ holds}
$$1=\sum_{i=0}^k(-1)^i{k\choose i}2^{k-i}.$$

\noindent\textbf{Note 2.1.} \textit{The preceding equation is a
trivial consequence of Binomial theorem. But here  it is obtained
independently of this theorem.}

The preceding Proposition shows that counting  i-columns over all
subsets  of $X$ always produce the same result.

We shall now make some restrictions on the number of subsets of
$X$. Take $0\leq m_1\leq m_2\leq n$ fixed,  and consider the
submatrix $C$ of $A$ consisting of rows whose indices belong to
$I_0,$ and columns corresponding to those subsets of $X$ that have
$m,\;(m_1\leq m\leq m_2)$ elements.

Let $c$ be an i-column of $C.$ Define $I_0'=\{i\in
I_0:c_i=1\},\;|I_0'|=l.$

The number $\nu(c)$ is equal to the number of sets that contain
$\{x_i:i\in I_0'\},$ and do not intersect the sets $\{x_i:i\in
I_0\setminus I_0'\}.$ We thus have
$$\nu=\sum_{i=m_1-|I_0'|}^{m_2-|I_0'|}{n-|I_0|\choose i}.$$

On the other hand, for $I\subseteq I_0$ the number $C(I)$
corresponds to the number of sets that contain
 $\{x_i:i\in I\setminus  I_0'\},$ and do not intersect $\{x_i:i\in I\cap
I_0'\}.$ Its number is equal
$$\sum_{i_3=m_1-|I\setminus I_0'|}^{m_2-|I\setminus I_0'|}{n-|I|\choose
i_3}.$$

It follows  that the formula (\ref{rv}) may be applied. We thus
have

 \noindent\textbf{Proposition 2.2.} \textit{ For $0\leq m_1\leq
m_2\leq n,$ and $0\leq l\leq k$ holds}

\begin{equation}\label{eq1}\sum_{i=m_1-l}^{m_2-l}{n-k\choose i}=\sum_{i_1=0}^{l}\sum_{i_2=0}^{k-l}\sum_{i_3=m_1-i_2}^{m_2-i_2}(-1)^{i_1+i_2}{l\choose
i_1}{k-l\choose i_2} {n-i_1-i_2\choose i_3}\end{equation}

In the special case when one takes $k=l,\;m_1=m_2=m$ we obtain

\noindent\textbf{Corollary 2.1.} \textit{ For arbitrary
nonnegative integers $m,n,k$ holds}
\begin{equation}\label{eq2}{n-k\choose m-k}=\sum_{i=0}^{k}(-1)^{i}{k\choose
i} {n-i\choose m}.\end{equation}

\noindent\textbf{Note 2.2.} \textit{The preceding is one of the
best known binomial identities. It appears in the book $[1]$ in
many different forms.}

Taking $m_1=m_2=m,$ in (\ref{eq1}) one gets

\noindent\textbf{Corollary 2.2.} \textit{ For arbitrary
nonnegative integer $m,n,k,l,\;(l\leq k)$ holds}
\begin{equation}\label{eq3}{n-k\choose m-l}=\sum_{i_1=0}^{l}\sum_{i_2=0}^{k-l}(-1)^{i_1+i_2}{l\choose
i_1}{k-l\choose i_2} {n-i_1-i_2\choose m-i_2},\end{equation}

For $l=0$ we obtain

\begin{equation}\label{eq4}{n-k\choose m}=\sum_{i=0}^{k}(-1)^{i}{k\choose i} {n-i\choose m-i},\end{equation}
which is only  another form of (\ref{eq2}).

Taking $n=2k,\;l=k$  in (\ref{eq3})we obtain
$${k\choose
m-k}=\sum_{i_1=0}^{k}(-1)^{i_1}{k\choose i_1}{2k-i_1\choose m}.$$
Substituting $k-i_1$ by $i$  we obtain

\noindent\textbf{Corollary 2.3. Klee's identity},([2],p.13)
$$(-1)^k{k\choose
m-k}=\sum_{i=0}^{k}(-1)^{i}{k\choose i}{k+i\choose m}.$$

From (\ref{eq1}) we may obtain different formulae for partial sums
of binomial coefficients, that is, for the numbers of Bernoulli's
triangle. For instance, taking $l=0,\;m_1=0,\;m_2=m$ we obtain

\noindent\textbf{Corollary 2.4.} For any $0\leq m\leq n$ and
arbitrary nonnegative integer $k$ holds
\begin{equation}\label{eq4}\sum_{i=0}^{m}{n\choose i}=\sum_{i_1=0}^{k}\sum_{i_2=0}^{m-i_1}(-1)^{i_1}{k\choose i_1} {n+k-i_1\choose i_2}.\end{equation}

\noindent\textbf{Note 2.3.} \textit{The number $k$ in the
preceding equation may be considered as a free variable that
 takes nonnegative integer values. Specially, for $k=1$ the equation
represents the standard recursion formula for the numbers of
Bernoulli's triangle.}

Taking $k=l=m_1,\;m_2=m$ one obtains

\begin{equation}\label{eq5}\sum_{i=0}^{m}{n\choose i}=\sum_{i_1=0}^{k}\sum_{i_2=k}^{m+k}(-1)^{i_1}{k\choose
i_1} {n+k-i_1\choose i_2}\end{equation}

\noindent\textbf{Note 2.4.} \textit{ The formulae $(\ref{eq4})$
and $(\ref{eq5})$ differs in the range of the index $i_2.$}

\vspace{1cm}

\noindent\textbf{References}

\vspace{0.2cm}

\noindent[1] J. Riordan, Combinatorial Identities. New York:
Wiley, 1979.

\end{document}